\documentclass[12pt,reqno]{amsart}
\usepackage{amssymb,delarray}
\usepackage{amsfonts}
\usepackage{epsfig}
\usepackage[all]{xy}

\makeindex{}

\newtheorem{thm}{Theorem}%[section]
\newtheorem{lem}%[thm]
{Lemma}

\newtheorem{cor}%[thm]
{Corollary}

{\catcode`\@=11
\gdef\n@te#1#2{\leavevmode\vadjust{%
 {\setbox\z@\hbox to\z@{\strut#1}%
  \setbox\z@\hbox{\raise\dp\strutbox\box\z@}\ht\z@=\z@\dp\z@=\z@%
  #2\box\z@}}}
\gdef\leftnote#1{\n@te{\hss#1\quad}{}}
\gdef\rightnote#1{\n@te{\quad\kern-\leftskip#1\hss}{\moveright\hsize}}
\gdef\?{\FN@\qumark}
\gdef\qumark{\ifx\next"\DN@"##1"{\leftnote{\rm##1}}\else
 \DN@{\leftnote{\rm??}}\fi{\rm??}\next@}}

\begin{document}
\baselineskip=14.pt plus 2pt %13.7pt plus 2pt

\title[On Chisini's Conjecture. II]{On Chisini's Conjecture. II}
\author[Vik.S.~Kulikov]{Vik. S.~Kulikov}
\address{Steklov Mathematical Institute\\
Gubkina str., 8\\
119991 Moscow \\
Russia} \email{kulikov@mi.ras.ru}

%\curraddr{}
\dedicatory{} \subjclass{}
\thanks{The work  was partially supported
by the RFBR  (05-01-00455),   NWO-RFBR 047.011.2004.026
%(РФФИ 05-02-89000-НВО_а)
(RFBR 05-02-89000-$NWO_a$), INTAS (05-1000008-7805), and by
RUM1-2692-MO-05. }
%\date{This version: July 2006. }
\begin{abstract}
It is proved that if $S\subset \mathbb P^N$ is a smooth projective
surface and $f:S\to \mathbb P^2$ is a generic linear projection
branched over a cuspidal curve $B\subset \mathbb P^2$, then the
surface $S$ is determined uniquely up to an isomorphism of $S$ by
the curve $B$.
\end{abstract}

\maketitle
%\pagenumbering{roman}
\setcounter{tocdepth}{2}
%\setcounter{tocdepth}{1}
%\tableofcontentsSymplectic isotopy of sections

%\setcounter{section}{6}

\def\st{{\sf st}}

\setcounter{section}{1}
%\section{Introduction}
Let $B\subset \mathbb P^2$ be an irreducible plane  algebraic curve
over $\mathbb C$ with ordinary cusps and nodes, as the only
singularities. Denote by $2d$ the degree of $B$, and let $g$ be
the genus of its desingularization, $c%= \# \{ \mbox{cusps of}\, B\}
$ the number of its cusps, and $n%= \# \{ \mbox{nodes of} \, B\}
$ the number of its nodes. A curve $B$ is called {\it the
discriminant curve of a generic covering  of the projective plane}
if there exists a finite morphism $f:S\to \mathbb P^2$, $\deg f\geq
3$, satisfying the following conditions:

$(i)$ $S$ is a non-singular irreducible projective surface;

$(ii)$ $f$ is unramified over $\mathbb P^2 \setminus B$;

$(iii)$  $f^{*}(B)=2R+C$, where $R$ is a non-singular irreducible
reduced curve and a curve $C$ is reduced;

$(iv)$  $f_{\mid R}:R\to B$ coincides with the normalization of
$B$. \newline Such $f$ is called {\it a generic covering of the
projective plane} $\mathbb P^2$.

A generic covering $f:S\to \mathbb P^2$ is called a {\it generic
projection} if the surface $S$ is embedded in some projective
space $\mathbb P^N$ and $f=\text{pr}_{\mid S}$ is a restriction to
$S$ of a linear projection $\text{pr}:\mathbb P^N\to \mathbb P^2$.

Chisini's Conjecture (see \cite{Chi}) claims that if $f:S\to
\mathbb P^2$ is a generic covering of the projective plane of
$\deg f \geq 5$ then $f$ is determined uniquely up to an
isomorphism of $S$ by its discriminant curve.

It was proved in \cite{Ku} that Chisini's Conjecture holds for the
discriminant curve $B$ of a generic covering $f:S\to \mathbb P^2$
if
\begin{equation}
\deg f>\frac{4(3d+g-1)}{2(3d+g-1)-c}. \label{in}
\end{equation}
Furthermore, it was observed in \cite{Nem} that, by Bogomolov --
Miaoka -- Yau inequality, the right side of inequality ($1$) takes
the values less then $12$, that is, Chisini's Conjecture holds for
the discriminant curves of the generic coverings of degree greater
than $11$. Besides, also it was shown in \cite{Nem} that if $S$ is a
surface of non-general type, then Chisini's Conjecture holds for the
discriminant curves of the generic coverings $f:S\to\mathbb P^2$ if
$\deg f\geq 8$.

The aim of the article is to prove the following theorem.

\begin{thm} \label{thm} Let $f:S\to \mathbb P^2$ be a generic projection. Then
the generic covering $f$ is uniquely determined up to an
isomorphism of $S$ by its discriminant curve $B\subset \mathbb
P^2$ except the case when $S\simeq \mathbb P^2$ is embedded in
$\mathbb P^5$ by the polynomials of degree two {\rm (}the Veronese
embedding of
 $\mathbb P^2$ in $\mathbb P^5${\rm )} and $f$ is the restriction
to $S$ of a linear projection $\text{pr}:\mathbb P^5\to \mathbb
P^2$.
\end{thm}

\proof To prove Theorem, we will show that inequality (\ref{in})
does not hold only for the discriminant curves of two continuous
families of generic projections onto the projective plane, and after
that we show that for one of these exceptional families, the generic
coverings $f:S\to \mathbb P^2$ are uniquely determined by their
discriminant curves  and the generic projections of the second
exceptional family are the generic projections of $S\simeq \mathbb
P^2$ embedded in $\mathbb P^5$ by the Veronese embedding.

For this purpose, consider a generic projection $f:S\to \mathbb
P^2$, where $S$ is a non-singular surface embedded in $\mathbb P^N$.
Let $\deg S=m$ be the degree of the embedding $S\subset \mathbb P^N$
and $\text{pr}:\mathbb P^N\to \mathbb P^2$ be a linear projection
such that $f=\text{pr}_{\mid S}$. We have $\deg f=\deg S=m$.

Any linear projection $\mathbb P^N\to \mathbb P^2$ is determined by
its center $\mathbb P^{N-3}\subset \mathbb P^N$. Therefore the set
of linear projections $\mathbb P^N\to \mathbb P^2$ is parameterized
by the points of the Grassmanian $\text{Gr}(N-3,N)$. Let $u_0\in
\text{Gr}(N-3,N)$ be a point for which the generic covering
$f=\text{pr}_{u_0\mid S}$ is the restriction of the projection
$\text{pr}=\text{pr}_{u_0}$. There is a Zariski open subset $U_S$ of
the Grassmanian $\text{Gr}(N-3,N)$ such that for each $u\in U_S$ the
restriction $f_u$ of the corresponding linear projection
$\text{pr}_u$ to $S$ is a generic covering of the projective plane.
The set $U_S$ is non-empty, since, by assumption, $u_0\in U_S$. For
$u\in U_S$, the discriminant curves $B_u$ of the generic coverings
$f_u$ have the same genus $g$, the same degree $\deg B_u=2d$, and
the same numbers $c$ and $n$ of the cusps and nodes. Therefore
inequality (\ref{in}) either holds or does not hold simultaneously
for all $f_u$, $u\in U_S$, and, consequently, any point of $U_S$ can
be taken as the point $u_0$ in order to check inequality (\ref{in}).

By Theorem 3 in \cite{Mo}, there is a non-empty Zariski open subset
$V_S\subset \text{Gr}(N-4,N)$ such that  for each $v\in V_S$ the
image $\overline S=\text{pr}_v(S)$ of $S$ under the linear
projection $\text{pr}_v:\mathbb P^N\to \mathbb P^3$ has only
ordinary singular points (that is, singular points given locally by
one of the following equations: $xy=0$ (a double curve), $xyz=0$ (a
triple point), and $x^2=y^2z$ (a pinch)).

Consider the flag manifold $F=F(N-4,N-3,N)$ of the linear subspaces
$\mathbb P^{N-4}\subset\mathbb P^{N-3}$ in $\mathbb P^N$. We have
two natural projections $p_1:F\to \text{Gr}(N-3,N)$ and $p_2:F\to
\text{Gr}(N-4,N)$. Obviously, the intersection
$W_S=p_1^{-1}(U_S)\cap p_2^{-1}(V_S)$ of two non-empty Zariski open
subsets $p_1^{-1}(U_S)$ and $p_2^{-1}(V_S)$ is a non-empty Zariski
open subset of $F$. Therefore, without loss of generality, we can
assume that the generic covering $f$ coincides with $f_u$ for some
$u\in U_S$ for which there is $w\in W_S$ such that $p_1(w)=u$, that
is, $\text{pr}_u$ can be decomposed into the composition of two
projections: the projection $\text{pr}_{p_2(w)}$ and a projection
$\overline{\text{pr}}:\mathbb P^3\to \mathbb P^2$ such that
$\overline S=\text{pr}_{p_2(w)}(S)$ is a surface in $\mathbb P^3$ of
degree $\deg \overline S=\deg S$ with ordinary singular points.
Denote by $f_1:S\to \overline S$ the restriction of
$\text{pr}_{p_2(w)}$ to $S$ and by $f_2:\overline S\to \mathbb P^2$
the restriction of $\overline{\text{pr}}$ to $\overline S$. The
morphism  $f_1$ is birational. We have $f=f_2\circ f_1$.

Denote by $D\subset \overline S$ the double curve of $\overline S$,
$D=D_1\cup \dots\cup D_u$, where $D_i$, $i=1,\dots,u$, are the
irreducible components of $D$. Let $g_i$ and $d_i$ be respectively
the genus and the degree of the curve $D_i$. Put
$\displaystyle\overline g=\sum_{i=1}^u g_i$ and
$\displaystyle\overline d=\sum_{i=1}^u d_i$. Denote by $t$ the
number of triple points of $\overline S$. Note that $0\leq u\leq
\overline d$ and $\overline g\geq 0$.

We have (see, for example, \cite{G-H})
\begin{equation} \label{K2}
K^2_S=m(m-4)^2-(5m-24)\overline d-4(u-\overline g)+9t,
\end{equation}
\begin{equation} \label{E}
e(S)=m^2(m-4)+6m-(7m-24)\overline d-8(u-\overline g)+15t,
\end{equation}
where $K_S$ is the canonical class of $S$ and $e(S)$ is its
topological Euler characteristic.  On the other hand, since $\deg
f=\deg S=m$ for a generic projection $f=\text{pr}_{\mid S}$, we have
(see Lemmas 6 and 7 in \cite{Ku})
\begin{equation} \label{k2}
K^2_S=9m-9d+g-1,
\end{equation}
\begin{equation} \label{e}
e(S)=3m+2(g-1)-c.
\end{equation}

\begin{lem} \label{dbard} We have
\begin{equation} \label{d} 2d=m(m-1)-2\overline d, \end{equation}
\begin{equation} \label{ind}
\overline d\leq \frac{(m-1)(m-2)}{2}.\end{equation}
\end{lem}

\proof Let $L$ be a generic line in $\mathbb P^2$ and $\overline
L=f_2^{-1}(L)$ its preimage. Then $\overline L$ is an irreducible
plane curve of degree $m$ having $\overline d$ nodes as its
singular points. Therefore its genus $g(\overline L)$ is equal to
$\frac{(m-1)(m-2)}{2}-\overline d$ and inequality (\ref{ind})
follows from the inequality $g(\overline L)\geq 0$.

The covering $f_{2\mid \overline L}:\overline L\to L$ is a morphism
of degree $m$ and it is branched at $2d=(L,B)_{\mathbb P^2}=\deg B$
points. Therefore, by Hurwitz formula, $2g(\overline L)-2=-2m +2d$.
Thus, we have
$$-2m+2d=(m-1)(m-2)-2\overline d-2,$$ that is, $2d=m(m-1)-2\overline d$.
\qed \\

It follows from equalities (\ref{K2}) -- (\ref{d}) that
\begin{equation} \label{g}
g-1=m\frac{2m^2-7m+5}{2}-5(m-3)\overline d-4(u-\overline g)+9t,
\end{equation}
\begin{equation} \label{c}
c=m(m-1)(m-2)-3(m-2)\overline d+3t. \end{equation}

Substituting equalities (\ref{d}), (\ref{g}), and (\ref{c}) in
inequality (\ref{in}) and performing evident transformations, it is
easy to show that inequality (\ref{in}) is equivalent to the
following inequality
\begin{equation} \label{in2}
(m-2)[m(m-1)(m-2)-(7m-24)\overline d-8(u-\overline g)]+3(5m-12)t>
0.
\end{equation}

Therefore, by  Theorem 1 in \cite{Ku}, to prove Theorem \ref{thm},
it suffices to show that if the inequality
\begin{equation} \label{nin2} (m-2)[m(m-1)(m-2)-(7m-24)\overline
d-8(u-\overline g)]+3(5m-12)t\leq 0
\end{equation}
holds for a surface $\overline S\subset \mathbb P^3$ with ordinary
singular points, then either $f:S\to\mathbb P^2$ is a projection
of the projective plane embedded in $\mathbb P^5$ by the Veronese
embedding or $f$ is uniquely determined up to an isomorphism of
$S$ by its discriminant curve $B$.

By the main result in \cite{Nem}, we can assume that $m\leq 11$.

\begin{lem} \label{gen}  Chisini's Conjecture holds for the discriminant curves of the
generic projections $f:S\to \mathbb P^2$ if $6\leq \deg S=m\leq
11$ and $K^2_S\leq 3e(S)$.
\end{lem}

\proof It follows from equalities (\ref{k2}), (\ref{e}), and the
inequality $K^2_S\leq 3e(S)$ that
\begin{equation} \label{BMYm}
3c\leq 9d+5(g-1),
\end{equation}
Assume that  Chisini's Conjecture does not hold for the discriminant
curve $B$ of a generic projection $f:S\to \mathbb P^2$, $\deg f=\deg
S=m$. Then the invariants of $B$ do not satisfy inequality
(\ref{in}), that is, these invariants satisfy the inequality
$$\frac{4(3d+g-1)}{2(3d+g-1)-c}\geq m,$$
or, equivalently,
\begin{equation} \label{nCh}
c\geq 2\frac{(m-2)}{m}(3d+g-1).
\end{equation}
It follows from inequalities (\ref{BMYm}) and (\ref{nCh}) that
$$6(m-2)[3d+(g-1)]\leq 3mc\leq m[9d+5(g-1)]$$
and hence
$$6(m-2)[3d+(g-1)]\leq m[9d+5(g-1)],$$
that is,
\begin{equation} \label{dgm}
g-1\geq \frac{9(m-4)}{12-m}d \end{equation} (by assumption, $m\leq
11$). Therefore, applying inequality (\ref{nCh}), we have
$$c\geq 2\frac{(m-2)}{m}(3d+g-1)\geq
2\frac{(m-2)}{m}(3d+\frac{9(m-4)}{12-m}d),$$ that is,
\begin{equation} \label{dcm}
c\geq \frac{12(m-2)}{12-m}d.
\end{equation}

Since $\frac{\deg B(\deg B-3)}{2}=c+n+g-1$ and $n\geq 0$, then
\begin{equation} \label{dgeq}
d(2d-3)\geq c+g-1.
\end{equation}
Therefore we have
$$d(2d-3)\geq c+g-1\geq \frac{12(m-2)}{12-m}d +
\frac{9(m-4)}{12-m}d$$ and hence
$$2d-3\geq \frac{12(m-2)}{12-m} + \frac{9(m-4)}{12-m}=\frac{21m-60}{12-m},$$ that is,
\begin{equation} \label{dgc}
d\geq  \frac{3(3m-4)}{12-m}.
\end{equation}

If $m=11$, then it follows from inequality (\ref{dgc}) that $d\geq
87$. On the other hand, by Lemma \ref{dbard}, $d\leq 55$.
Contradiction.

If $m=10$, then it follows from inequality (\ref{dgc}) that $d\geq
39$. Therefore, by Lemma \ref{dbard}, we have $\overline d\leq 6$.

On the other hand, inequality (\ref{nin2}) implies the inequality
$$8[720-46\overline d -8(u-\overline g)]+114t\leq 0$$
and, consequently, since $t\geq 0$, we should have
$$720- 46\overline d
-8(u-\overline g)\leq 0.$$ Therefore
$$720\leq 46\overline d +8(u-\overline g)\leq 54\overline d ,$$
since $u-\overline g\leq \overline d$. Finally, we obtain the
inequality $\overline d\geq \frac{720}{54}$, which contradicts the
inequality $\overline d\leq 6$.

If $m=9$, then it follows from inequalities (\ref{dgm}) and
(\ref{dgc}) that $d\geq 23$ and $g-1\geq 15d$. Therefore, by Lemma
\ref{dbard}, we have
\begin{equation} \label{god9}
g-1\geq 15(36-\overline d),
\end{equation}
\begin{equation} \label{od9}
\overline d\leq 28-23=5.
\end{equation}

It follows from inequality (\ref{nin2}) that
$$7[504-39\overline d -8(u-\overline g)]+99t\leq 0,$$
or, equivalently,
\begin{equation} \label{In9}
99t\leq 273\overline d +56(u-\overline g)-3528.\end{equation}

Equality (\ref{g}), in which we substitute $m=9$, and inequality
(\ref{god9}) imply the following inequality
$$468 -30 \overline d-4(u-\overline g)+9t\geq 15(36-\overline d), $$
or, equivalently,
\begin{equation} \label{og9}
9t\geq 15 \overline d +4(u-\overline g)+72.
\end{equation}
It follows from inequalities (\ref{In9}) and (\ref{og9}) that
$$[273\overline d +56(u-\overline g)-3528]\geq 11[15 \overline d +4(u-\overline g)+72],$$
that is, $4320\leq 108\overline d+12(u-\overline g)\leq 120\overline
d$, since $\overline g\geq 0$ and $u\leq \overline d$. Therefore
$\overline d\geq 36$, but this inequality contradicts inequality
(\ref{od9}).

If $m=8$, then it follows from inequality (\ref{dgm}) that
$g-1\geq 9d$. Therefore, by Lemma \ref{dbard}, we have
\begin{equation} \label{god8}
g-1\geq 9(28-\overline d),
\end{equation}
where $\overline d\leq 21$.

It follows from inequality (\ref{nin2}) that
$$6[336-32\overline d -8(u-\overline g)]+84t\leq 0,$$
or, equivalently,
\begin{equation} \label{In8}
7t\leq 16\overline d +4(u-\overline g)-168.\end{equation}

Equality (\ref{g}), in which we substitute $m=8$, and inequality
(\ref{god8}) imply the following inequality
$$308 -25 \overline d-4(u-\overline g)+9t\geq 9(28-\overline d), $$
or, equivalently,
\begin{equation} \label{og8}
9t\geq 16 \overline d +4(u-\overline g)-56.
\end{equation}
It follows from inequalities (\ref{In8}) and (\ref{og8}) that
$$7[16 \overline d +4(u-\overline g)-56]\leq  9[16\overline d +4(u-\overline g)-168],$$
that is, $1120\leq 32\overline d+8(u-\overline g)\leq 40\overline
d$, since $\overline g\geq 0$ and $u\leq \overline d$. Therefore
$\overline d\geq 28$, which contradicts the inequality $\overline
d\leq 21$.

If $m=7$, then it follows from inequality (\ref{dgm}) that
$g-1\geq \frac{27}{5}d$. Hence we have \begin{equation}
\label{gd7} g-1\geq \frac{27}{5}(21-\overline d),\end{equation}
since, by Lemma \ref{dbard}, $d=21-\overline d$ and $\overline
d\leq 15$.

Inequality (\ref{nin2}) can be written in the following form
\begin{equation} \label{In7}
69t\leq 125\overline d +40(u-\overline g)-1050.
\end{equation}

Equality (\ref{g}), in which we substitute $m=7$, and inequality
(\ref{gd7}) imply the following inequality
$$189 -20 \overline d-4(u-\overline g)+9t\geq \frac{27}{5}(21-\overline d), $$
or, equivalently,
\begin{equation} \label{og7}
45t\geq 73 \overline d+20(u-\overline g)-378. \end{equation} It
follows from inequalities (\ref{In7}) and (\ref{og7}) that
$$15[125\overline d+40(u-\overline g)-1050]\geq 23[73 \overline d+20(u-\overline g)-378],$$
that is, $7056\leq 196\overline d+140(u-\overline g)\leq
336\overline d$, since $\overline g\geq 0$ and $u\leq \overline d$.
Therefore $\overline d\geq \frac{7056}{336}=21$, which contradicts
the inequality $\overline d\leq 15$.

If $m=6$, then it follows from inequality (\ref{dgm}) that $g-1\geq
3d$. Hence we have \begin{equation} \label{gd6} g-1\geq
3(15-\overline d),\end{equation} since, by Lemma \ref{dbard},
$d=15-\overline d$  and $\overline d\leq 10$.

Inequality (\ref{nin2}) can be written in the following form
\begin{equation} \label{In6}
27t\leq 36 \overline d+16(u-\overline g)-240.
\end{equation}

Equality (\ref{g}), in which we substitute $m=6$, and inequality
(\ref{gd6}) imply the following inequality
$$105 -15 \overline d-4(u-\overline g)+9t\geq 45-3\overline d, $$
or, equivalently (multiplying by $3$),
\begin{equation} \label{og6}
27t\geq 36\overline d+12(u-\overline g)-180.\end{equation} It
follows from inequalities (\ref{In6}) and (\ref{og6}) that
$$36 \overline d+16(u-\overline g)-240\geq 36\overline d+12(u-\overline g)-180,$$
that is, $u-\overline g\geq 15$. On the other hand, we have
$u-\overline g\leq 10$, since $\overline g\geq 0$ and $u\leq
\overline d\leq 10$. Contradiction. \qed \\

By Theorem 2 in \cite{Nem}, if $\deg f\geq 8$ and $S$ is a surface
of non-general type, then Chisini's Conjecture holds for the
discriminant curve $B$ of any generic covering $f:S\to\mathbb P^2$.
It is well known (see the classification of the algebraic surfaces),
that if Bogomolov -- Miaoka -- Yau inequality does not take place
for an algebraic surface $S$, then $S$ is an irregular ruled surface
and in this case we have $K^2_S\leq 2e(S)$ and $K^2_S\leq -2$.
Therefore, by Lemma \ref{gen}, to prove Theorem \ref{thm}, it
suffices to consider only the following cases: $3\leq m\leq 7$ and
if $m=6$ or $7$, then $K_S^2\leq 2e(S)$ and $K^2_S\leq -2$.

Again, we assume that the invariants of the surface $\overline S$
satisfy inequality (\ref{nin2}).

{\it Case $m=3$.} In this case inequality (\ref{nin2}) has the
form
$$ 6+3\overline d -8(u-\overline g) +9t\leq 0.$$
It follows from inequality (\ref{ind}) that $\overline d\leq 1$
and we have two possibilities: either $\overline d=0$ and,
consequently, $u=\overline g=t=0$, or $\overline d=1$ and,
consequently, $u=1$, $\overline g=t=0$, since in this case $D$ is
a line in $\mathbb P^3$. It is easy to see that inequality
(\ref{nin2}) does not hold in both cases.

{\it Case $m=4$.} In this case inequality (\ref{nin2}) has the
form
$$ 2(24- 4\overline d -8(u-\overline g)) +24t\leq 0.$$
It follows from inequality (\ref{ind}) that $\overline d\leq 3$ and
we have three possibilities: either $\overline d\leq 2$ and,
consequently, $u\leq \overline d\leq 2$, $\overline g=t=0$, or
$\overline d=3$, $u=3$, $\overline g=0$, $t=1$, or $\overline d=3$,
$u=1$, $t=0$, and $\overline g=1$ or $0$. It is easy to see that
inequality (\ref{nin2}) holds only in the following two cases:
$u=\overline d=2$,  $\overline g=t=0$ and  $u=\overline d=3$,
$\overline g=0$, $t=1$. These exceptional cases will be investigated
at the end of the proof of Theorem \ref{thm}.

{\it Case $m=5$.}  Inequality (\ref{nin2}) has the following form
$$3[60-11\overline d-8(u-\overline g)]+39t\leq 0$$
or, equivalently,
\begin{equation}
\label{in5} 60+2t\leq 11(\overline d-t)+8(u-\overline
g).\end{equation}

By Theorem 11 in \cite{Ku}, Chisini's Conjecture holds for the
cuspidal curves $B$ of genus $g\leq 3$. Therefore, by equality
(\ref{g}), we should have $$g-1=50-10\overline d-4(u-\overline
g)+9t\geq 3$$ or, equivalently,
\begin{equation} \label{g5}
47-t\geq 10(\overline d-t)+4(u-\overline g).
\end{equation}

By Lemma \ref{dbard}, we have $u\leq \overline d\leq 6$. Therefore
$u-\overline g\leq 6$ and it follows from inequality (\ref{in5})
that
$$12+2t\leq 11(\overline d-t),$$
that is, \begin{equation} \label{a} \overline d-t\geq 2.
\end{equation}

Similarly, since $\overline d-t\leq 6$, it follows from inequality
(\ref{in5}) that
$$-6+2t\leq 8(u-\overline g)$$ and hence $u-\overline g\geq 0$.
Applying inequality (\ref{g5}), we have $47-t\geq 10(\overline
d-t)$, that is, $\overline d-t\leq 4$. Therefore, by inequality
(\ref{in5}), we obtain that $16+2t\leq 8(u-\overline g)$, that is,
$u-\overline g\geq 2$. Then, by inequality (\ref{g5}), we have
$39-t\geq 10(\overline d-t)$ and hence \begin{equation} \label{b}
\overline d-t\leq 3.
\end{equation} Now, it follows from inequality (\ref{in5}) that
$$27+2t\leq 8(u-\overline g)$$ and thus $u-\overline g\geq 4$.
Therefore $u\geq 4$ and hence $\overline d\geq 4$.

By inequalities (\ref{a}) and (\ref{b}), we have
$$ 2\leq \overline d-t\leq 3.$$
Let us consider the case $\overline d-t= 3$. It follows from
inequality (\ref{g5}) that \begin{equation} \label{t} 17-t\geq
4(u-\overline g). \end{equation} Therefore $u-\overline g\leq 4$.
Hence $u-\overline g=4$, $u=4$, and $\overline g =0$, since the
genera of the irreducible components of a curve of degree
$\overline d\leq 6$, having more then four irreducible components,
should be equal to zero. In addition, it follows from inequality
(\ref{t}) that $t\leq 1$. Therefore $t=1$ and $\overline d= 4$,
since $\overline d-t= 3$ and $\overline d\geq 4$. In this case, by
formulae (\ref{g}), (\ref{c}), and Lemma \ref{dbard}, the curve
$B$ should have the following invariants: $$\deg B =2d=12, \,\,
g=4, \,\, c=27,\,\, n=(2d-1)(d-1)-g-c=24.$$ But, it is impossible,
since in this case, by Pl\"{u}cker's formula, the dual curve
$\check{B}$ has degree $2d(2d-1)-3c-2n=3$ and therefore $\deg B$
should be less or equal
$$\deg \check{B}(\deg \check{B}-1)=3\cdot 2=6.$$

Let us consider the case $\overline d-t= 2$. It follows from
inequality (\ref{in5}) that
\begin{equation} \label{r} 38+2t\leq
8(u-\overline g). \end{equation}
Therefore $u-\overline g\geq 5$.
Hence $u\geq 5$ and $\overline g =0$. Now, by inequality
(\ref{r}), we should have $u=6$, since it follows from the
equality $\overline d-t= 2$ and the inequalities $6\geq \overline
d\geq u\geq 5$ that $t\geq 3$ and hence $38+2t\geq 44$. Therefore
we have only the following possibility:
$$u=6, \,\, \overline g=0, \, \, \overline d=6, \, \, t=4.$$ But,
these values of $u$, $\overline g$, $\overline d$, and $t$ do not
satisfy inequality (\ref{g5}).

{\it Case $m=6$ and $K^2_S\leq 2 e(S)$.} Applying formulae
(\ref{k2}) and (\ref{e}), we obtain the inequality
\begin{equation} \label{nBMY6} 2c\leq
9d+3(g-1)-18.
\end{equation}

Inequality (\ref{nCh}) can be written in the following form
\begin{equation} \label{in6}
3c\geq 4(3d+g-1).
\end{equation}
It follows from inequalities (\ref{in6}) and (\ref{nBMY6}) that
$$24d+8(g-1)\leq 6c\leq 27d+9(g-1)-54,$$
that is, $3d+g-1\geq 54$ and, since $d=15-\overline d$, we have
\begin{equation} \label{6g}
g-1\geq 54-3(15-\overline d)=9+3\overline d.
\end{equation}

By assumption, the invariants of $\overline S$ should satisfy
inequality (\ref{nin2}), where $m=6$, that is, they satisfy
inequality (\ref{In6}).

Equality (\ref{g}), in which we substitute $m=6$, and inequality
(\ref{6g}) imply the following inequality
$$105 -15 \overline d-4(u-\overline g)+9t\geq 9+3\overline d, $$
or, equivalently,
$$27t\geq 72\overline d+12(u-\overline g)-288.$$
By inequality (\ref{In6}), we have
$$36\overline d+16(u-\overline g)-240\geq 27t.$$
Therefore
$$36\overline d+16(u-\overline g)-240\geq 72\overline d+12(u-\overline g)-288,$$
that is, $12\geq 9\overline d-(u-\overline g)$. But, $u-\overline
g\leq \overline d$. Therefore $9\overline d-(u-\overline g)\geq
8\overline d$ and hence $3\geq 2\overline d$, that is, we should
have $\overline d\leq 1$, since $\overline d$ is an integer.

On the over hand, by inequality (\ref{In6}), we have
$$240\leq 9\overline d+27(\overline d-t)+16(u-\overline g)\leq
52\overline d,$$ since $\overline d-t\leq \overline d$ and
$u-\overline g\leq \overline d$. Therefore we should have
$\overline d\geq 5$. Contradiction.

{\it Case $m=7$ and $K^2_S\leq 2e(S)$,  $K^2_S\leq -2$.} Applying
formulae (\ref{k2}), (\ref{e}), and inequality  $K^2_S\leq 2e(S)$,
we obtain the inequality
\begin{equation} \label{nBMY7} 2c\leq
9d+3(g-1)-21.
\end{equation}

We have $K^2_S\leq -2$. Therefore it follows from formula (\ref{K2})
that
$$ K^2_S=7\cdot 9 - 25\overline d-4(u-\overline
g)+9t\leq -2.
$$
Hence
$$65\leq 65+9t\leq 25\overline d+4(u-\overline
g)\leq 29 \overline d, $$ since $t\geq 0$ and $u-\overline g\leq
\overline d$. Thus, we have
\begin{equation} \label{7d3}
\overline d\geq 3
\end{equation}

Inequality (\ref{nCh}) can be written the following form
\begin{equation}
\label{in7} 7c\geq 10(3d+g-1).
\end{equation} It  follows from inequalities (\ref{in7}) and (\ref{nBMY7}) that
$$60d+20(g-1)\leq 14c\leq 63d+21(g-1)-147,$$
that is,
\begin{equation} \label{gdc} 3d+g-1\geq 147
\end{equation}
and, since $d=21-\overline d$, we have
\begin{equation}
\label{7g1} g-1\geq 147-3(21-\overline d)=84+3\overline d.
\end{equation}
Therefore inequality (\ref{7d3}) implies the following inequality
\begin{equation} \label{7g}
g-1\geq 93.
\end{equation}

It follows from inequalities (\ref{in7}) and (\ref{gdc}) that
$c\geq 210$ and, by inequality (\ref{dgeq}), we have
$$d(2d-3)\geq c+g-1\geq 210+93=303$$ and hence
$d\geq\frac{3+\sqrt{2433}}{4}>13$, that is,
\begin{equation} \label{7d7} d\geq
14,\end{equation} since $d$ is an integer. Therefore
\begin{equation} \label{od7}
\overline d=21-d\leq 7.
\end{equation}

By assumption, the invariants of $\overline S$ should satisfy
inequality (\ref{nin2}), where $m=7$, that is, they satisfy
inequality (\ref{In7}). It follows from inequality (\ref{In7})
that
$$210-25\overline d -8(u-\overline g)\leq 0,$$
since $t\geq 0$. Therefore
$$210\leq 25\overline d +8(u-\overline g)\leq 33\overline d,$$
since $u-\overline g\leq \overline d$, and hence $\overline d\geq
\frac{210}{33}=6+\frac{4}{11}$, that is, $\overline d\geq 7$,
since $\overline d$ is an integer. Applying inequality
(\ref{od7}), we should have $\overline d=7$.

Equality (\ref{g}), in which we substitute $m=7$, and inequality
(\ref{7g1}) imply the following inequality
$$181 -20 \overline d-4(u-\overline g)+9t\geq 84+3\overline d .$$
Therefore \begin{equation} \label{xx} 9t\geq 64+4(u-\overline
g),\end{equation} since $\overline d=7$, and by (\ref{In7}), we
have $$125\overline d+40(u-\overline g)-1050\geq 69t,$$ that is,
\begin{equation} \label{xxx}
40(u-\overline g)-175\geq 69t. \end{equation}
Combining
inequalities (\ref{xx}) and (\ref{xxx}), we obtain the inequality
$$3[40(u-\overline g)-175]\geq 23[64+4(u-\overline g)],$$ that is,
$28(u-\overline g)\geq 3\cdot 175+23\cdot 64=1997$. But, on the
other hand, $u-\overline g\leq \overline d=7$. Contradiction.

Let us return to the last two non-investigated cases when $m=4$
and $u=\overline d=2$,  $\overline g=t=0$, or   $u=\overline d=3$,
$\overline g=0$, $t=1$.

Consider the case $m=4$ and $u=\overline d=2$,  $\overline g=t=0$.
By formulae (\ref{g}), (\ref{c}), and (\ref{d}), we have $d=4$,
$g=1$, and $c=12$. Therefore the number $n$ of nodes of $B$ is
equal to $d(2d-3)-c-g+1=8$.

Suppose that there is another generic covering $f_2:S_2\to \mathbb
P^2$ with the same discriminant curve $B$ which is non-equivalent
to the generic projection $f$. Then, by Theorem 1 in \cite{Ku},
$$\deg f_2\leq \frac{4(3d+g-1)}{2(3d+g-1)-c}=4.$$ Since $S_2$ is a
non-singular surface and the discriminant curve $B$ of $f_2$ has
nodes, $\deg f_2$ can not be equal to $3$ and hence $\deg f_2=4$.

Put $S_1=S$, $R_1=R$, $C_1=C$, and $f^*_2(B)=2R_2+C_2$, where
$R_2$ is the ramification locus of $f_2$.

Consider the fibred product
$$S_1\times _{\mathbb P^2}S_2=\{ \,
(x,y)\in S_1\times S_2\, \,  \mid \, \, f_1(x)=f_2(y)\, \, \} $$ and
let $X=\widetilde{S_1\times _{\mathbb P^2}S_2}$ be the normalization
of $S_1\times _{\mathbb P^2}S_2$. Denote  the corresponding natural
morphisms by $g_{1}:X\to S_1$, $g_{2}:X\to S_2$, and $f_{1,2}:X\to
\mathbb P^2$. We have $\deg g_1=\deg f_2=4$, $\deg g_2=\deg f_1=4$,
and $\deg f_{1,2}=\deg g_1\cdot\deg f_1=16$. By Propositions 2 and 3
in \cite{Ku},  $X$ is an irreducible non-singular surface.

Let $\widetilde R \subset X$ be a curve $g_1^{-1}(R_1)\cap
g_2^{-1}(R_2)$, $\widetilde C=g_1^{-1}(C_1)\cap g_2^{-1}(C_2)$,
$\widetilde C_1=g_1^{-1}(R_1)\cap g_2^{-1}(C_2)$, and $\widetilde
C_2=g_1^{-1}(C_1)\cap g_2^{-1}(R_2)$.

By Proposition 4 in \cite{Ku}, we have
$$
\begin{array}{cll}
\widetilde R^2 & = & 2(3d+g-1)-c=12 \, ,  \\
 \widetilde C_1^2 & = & (\deg f_1-2)(3d+g-1)-c =12  \, , \\
\widetilde C_2^2 & = & (\deg f_2-2)(3d+g-1)-c =12  \, , \\
(\widetilde R,\widetilde C_i) & = & c=12\, \hspace{1cm}
\mbox{for}\, \, i=1,\, 2
\end{array}
$$
and applying the same arguments which was used in the proof of
Proposition 4 in \cite{Ku}, it can be easily shown that the
intersection number $$(\widetilde C_1,\widetilde C_2)=c+2n=28.$$
Therefore the determinant
$$
\left|
\begin{array}{cc}
\widetilde R^2 & (\widetilde R,\widetilde C_1) \\
(\widetilde C_1,\widetilde R) &   \widetilde C_1^2
\end{array}
\right| =0
$$
and, consequently, by Hodge's Index Theorem, the classes
$[\widetilde C_1]$ and $[\widetilde R]$ of the curves $\widetilde
C_1$ and $\widetilde R$ in the Neron -- Severi group
$\text{NS}(X)$ of the divisors on $X$ up to numerical equivalence
are linear dependent. Since $\widetilde R^2=\widetilde C_1^2$,
then $[\widetilde R]=[\widetilde C_1]$ in $\text{NS}(X)$. Applying
the same arguments, we obtain that $[\widetilde R]=[\widetilde
C_2]$ and hence $[\widetilde C_2]=[\widetilde C_1]$ in
$\text{NS}(X)$. Therefore the intersection number $(\widetilde
C_2,\widetilde C_1)$ must be equal to $\widetilde C_1^2=12$. On
the other hand, $(\widetilde C_2,\widetilde C_1)=28$.
Contradiction.

To complete the proof of Theorem \ref{thm}, note that the last
case when $m=4$, $u=\overline d=3$, $\overline g=0$, $t=1$
corresponds to a generic projection $f:S\to \mathbb P^2$, where
$S\simeq \mathbb P^2$ is embedded in $\mathbb P^5$ by the
polynomials of degree two {\rm (}the Veronese embedding of
 $\mathbb P^2$ in $\mathbb P^5${\rm )} and $f$ is the restriction
to $S$ of a linear projection $\text{pr}:\mathbb P^5\to \mathbb
P^2$ (see, for example, \cite{G-H}). In this case, $B\subset
\mathbb P^2$ is the dual curve of a smooth cubic, $\deg B=6$,
$c=9$, and $B$ is the discriminant curve of four non-equivalent
generic coverings of $\mathbb P^2$ (\cite{Chi}, \cite{Cat}). Three
of them have degree four, and the last one has degree three. \qed

\begin{cor} Let $S_i$ be non-singular surfaces, $i=1,2$, and $S_i\subset \mathbb P^{N_i}$ two
embeddings given by complete linear systems of divisors on $S_i$.
Suppose that these embeddings do not coincide with the Veronese
embedding of $\mathbb P^2$ in $\mathbb P^5$. Let
$f_i=\text{pr}_{i\mid S_i}:S_i\to \mathbb P^2$ be two generic
coverings ramified over the same cuspidal curve $B$, where
$\text{pr}_i:\mathbb P^{N_i}\to \mathbb P^2$ are linear projections.
Then $N_1=N_2=N$ and there is a linear transformation $h:\mathbb
P^N\to \mathbb P^N$ such that $h(S_1)=S_2$ and $f_1=f_2\circ h$.
\end{cor}

\proof  Denote by $\overline L_i=f_i^{-1}(L)\subset S_i$, $i=1,2$,
the proper transform of a line $L$ in $\mathbb P^2$. By Theorem
\ref{thm}, there is an isomorphism $h:S_1\to S_2$ such that
$f_1=h\circ f_2$. Therefore $h(\overline L_1)=\overline L_2$ and
hence $h^*(\mathcal O_{S_2}(\overline L_2))=\mathcal
O_{S_1}(\overline L_1)$. Consequently, we have $$ N_1=\dim
H^0(S_1,\mathcal O_{S_1}(\overline L_1))=\dim H^0(S_1,\mathcal
O_{S_2}(\overline L_2))=N_2$$ and the isomorphism $h$ can be
defined by a linear transformation $\mathbb P^{N_1}\to \mathbb
P^{N_2}$ induced by $h^*:H^0(S_2,\mathcal O_{S_2}(\overline
L_2))\to H^0(S_1,\mathcal O_{S_1}(\overline L_1))$. \qed

Note also that if $f:S\to\mathbb P^2$ is a generic covering, $\deg
f=4$, branched over a cuspidal curve $B\subset \mathbb P^2$, $\deg
B=6$, $c=9$, then, by (\ref{k2}) and (\ref{e}), we have $K^2_S=9$
and $e(S)=3$. By Hurwitz formula, the genus of $f^{-1}(L)$, where
$L$ is a line in $\mathbb P^2$, is equal to $\frac{-2\deg
f+(L,B)}{2}+1=0$. Therefore $S\simeq \mathbb P^2$ and $f$ is given
by polynomials of degree 2. Hence, in the exceptional case of a
cuspidal curve $B\subset \mathbb P^2$, $\deg B=6$, $c=9$, each of
three non-equivalent generic coverings $f_i$, $\deg f_i=4$, ramified
over $B$, are generic projections of $\mathbb P^2$ embedded in
$\mathbb P^5$ by Veronese embeddings. It is easy to see that the
fourth exceptional generic covering $f_4:S\to\mathbb P^2$, $\deg
f_4=3$, is not a generic projection (see Case $m=3$ in the proof of
Theorem \ref{thm}). Therefore we have the following
\begin{cor} Let $f:S\to \mathbb P^2$ be a generic linear projection
branched over a cuspidal curve $B\subset \mathbb P^2$.
 Then $S$ is determined uniquely (up to isomorphism) by the curve $B$.
\end{cor}

\ifx\undefined\bysame
\newcommand{\bysame}{\leavevmode\hbox to3em{\hrulefill}\,}
\fi

\end{document}